\documentclass[11pt]{article}

\usepackage[T1]{fontenc}
\usepackage[utf8]{inputenc}
\usepackage{lmodern}
\usepackage[margin=1.1in]{geometry}
\usepackage{amsmath,amssymb,amsthm}
\usepackage{booktabs}
\usepackage{array}
\usepackage{microtype}
\usepackage[colorlinks=true,linkcolor=blue,citecolor=blue,urlcolor=blue]{hyperref}

\newtheorem{theorem}{Theorem}[section]
\newtheorem{lemma}[theorem]{Lemma}
\newtheorem{proposition}[theorem]{Proposition}
\newtheorem{corollary}[theorem]{Corollary}
\theoremstyle{definition}

\theoremstyle{remark}
\newtheorem{remark}[theorem]{Remark}

\newcommand{\Ind}[1]{\mathbf{1}\!\left[#1\right]}
\newcommand{\Pmn}{P_m \times P_n}
\newcommand{\dis}{M}
\DeclareMathOperator{\dist}{d}

\title{The maximum length of a Hamiltonian path\\ in the grid graph $P_m \times P_n$}
\author{Jiajun Yuan\\[3pt] \small Independent researcher\\ \small \texttt{df11g0203@outlook.com}}
\date{}

\begin{document}

\maketitle

\begin{abstract}
\noindent
Label the $mn$ vertices of the grid graph $P_m \times P_n$ bijectively by
$1, \dots, mn$ and read the labels in increasing order; this produces a
Hamiltonian path whose length in the Manhattan metric is denoted
$\dis(\Pmn)$.  McNeil conjectured, on the basis of computer evidence, that
$\dis(P_n \times P_n) = n^3 - 3$ for even $n$ and $n^3 - n - 1$ for odd $n$.
Fried recently determined $\dis(\Pmn)$ up to an additive error of one.  We
determine it exactly, which settles McNeil's conjecture.  The proof rests on
the cut-counting bound for multisets.  Fried states that bound for
multiplicities at least two; we give a form valid for every multiplicity, with
a correction term that is invisible in the planar case but not in higher
dimensions.
\end{abstract}



\section{Introduction}

\subsection{The problem}

Let $m, n \ge 2$ and let $\Pmn$ be the Cartesian product of the paths
$P_m$ and $P_n$.  We identify its vertex set with
$V = \{(i,j) : i \in [m],\, j \in [n]\}$, where $[k] = \{1,\dots,k\}$, and we
write $N = mn$ and
\[
  \dist\bigl((i_1,j_1),(i_2,j_2)\bigr) = |i_1 - i_2| + |j_1 - j_2|
\]
for the Manhattan distance.

A \emph{labelling} is a bijection $\sigma \colon V \to [N]$.  Put
$v_t = \sigma^{-1}(t)$; the labelling induces the Hamiltonian path
$v_1, v_2, \dots, v_N$, and we are interested in its length
\[
  L(\sigma) \;=\; \sum_{t=1}^{N-1} \dist(v_t, v_{t+1}),
  \qquad
  \dis(\Pmn) \;=\; \max_{\sigma} L(\sigma).
\]
For a general graph $G$ the analogous quantity, the maximum length of a
Hamiltonian path in the metric of $G$, is the \emph{maximum travelling
salesman problem}; see \cite{GutinPunnen} for a survey.  (In the
\emph{maximum scatter} variant one maximises the shortest edge rather than
the total length; see \cite{Arkin}.)  For grid graphs in general, that is,
for induced subgraphs of the rectangular grid, even deciding whether a
Hamiltonian path exists is NP-complete \cite{IPS}; the full rectangles
$P_m \times P_n$ always admit one, and the question here is how long the
longest one can be.  In one dimension the problem is the study of
\emph{disorders} of permutations \cite{BGV}.

\subsection{Previous results}

Three results are known.

\begin{enumerate}
\item For $m = 1$ the grid degenerates to a path and the problem is that of
  maximising $\sum_{t} |\sigma_{t+1} - \sigma_t|$ over permutations of $[n]$.
  Bulteau, Giraudo and Vialette \cite[Theorem 8]{BGV} proved
  \begin{equation}\label{eq:oned}
    \dis(P_n) \;=\; \left\lfloor \frac{n^2}{2} \right\rfloor - 1 .
  \end{equation}

\item For $m = n$, McNeil computed the values for $n \le 19$
  (sequence A179094 of the OEIS \cite{OEIS}) and conjectured
  \begin{equation}\label{eq:mcneil}
    \dis(P_n \times P_n) \;=\;
    \begin{cases}
      n^3 - 3, & n \ \text{even},\\[2pt]
      n^3 - n - 1, & n \ \text{odd}.
    \end{cases}
  \end{equation}

\item Fried \cite[Theorem 1]{Fried}\footnote{References to \cite{Fried} are to the published version; the numbering in the arXiv preprint may differ.} proved that, for $m,n \ge 3$,
  \[
    \dis(\Pmn) \in \{r, r+1\},
  \]
  where
  \[
    r \;=\;
    \begin{cases}
      \dfrac{mn(m+n)}{2} - 3, & m,n \ \text{both even},\\[6pt]
      \dfrac{mn(m+n)}{2} - \dfrac{m+n}{2} - 1, & m,n \ \text{both odd},
    \end{cases}
  \]
  that $\dis(\Pmn) = \dfrac{mn(m+n)}{2} - \dfrac{n}{2} - 1$ when $m$ is odd
  and $n$ is even, and that $\dis(P_2 \times P_n) \in \{r, r+1\}$ with
  $r = (n+1)^2 - 4$ for every $n \ge 2$.
\end{enumerate}

Thus \cite{Fried} pins $\dis(\Pmn)$ down to within one in three of the four
families, and the remaining ambiguity is precisely the one that separates
\eqref{eq:mcneil} from its negation.

\subsection{Results of this note}

Our main result removes the ambiguity in all cases.

\begin{theorem}[Main theorem]\label{thm:main}
Let $m, n \ge 2$.  Then
\[
  \dis(\Pmn) \;=\;
  \begin{cases}
    \dfrac{mn(m+n)}{2} - 3, & m, n \ \text{both even},\\[8pt]
    \dfrac{mn(m+n)}{2} - \dfrac{m+n}{2} - 1, & m, n \ \text{both odd},\\[8pt]
    \dfrac{mn(m+n)}{2} - \dfrac{n}{2} - 1, & m \ \text{odd},\ n \ \text{even},\\[8pt]
    (n+1)^2 - 4, & \min(m,n) = 2 .
  \end{cases}
\]
\end{theorem}

In the third line $n$ is the even coordinate; by the symmetry
$\dis(\Pmn) = \dis(P_n \times P_m)$ the roles of $m$ and $n$ may be
interchanged.  In the fourth line $n = \max(m,n)$.  The four lines cover all
$m,n \ge 2$: the first three assume $m,n \ge 3$, the last one $\min(m,n)=2$.
On the overlap the first and fourth lines agree, since
$\frac{2n(2+n)}{2} - 3 = (n+1)^2 - 4$.

\begin{corollary}\label{cor:mcneil}
McNeil's conjecture \eqref{eq:mcneil} holds.
\end{corollary}

\begin{proof}
Take $m = n$ in Theorem~\ref{thm:main}.  If $n$ is even the first line gives
$n^3 - 3$; if $n \ge 3$ is odd the second line gives
$\frac{n^2 \cdot 2n}{2} - n - 1 = n^3 - n - 1$.  The case $n = 2$ is covered by
the fourth line and gives $9 - 4 = 5 = 2^3 - 3$.
\end{proof}

The new ingredients are of two kinds.  First, the upper bound of
\cite{Fried} is obtained from a bound for multisets (\cite[Lemma 2]{Fried}),
which is stated there for multiplicities at least two.  We give a proof of the
sharp form, valid for every multiplicity, and note that at multiplicity one a
correction term is needed --- a situation that arises in dimensions at least
three but not in the planar case.  Second, the gap of one is closed by a rigidity argument:
if the upper bound were attained, then each of the two coordinate sequences
would have to be extremal, and extremal sequences are forced to alternate
across the middle cut; two such alternations cannot coexist.  The three cases
of Theorem~\ref{thm:main} in which the gap has to be closed are treated
separately in Section~\ref{sec:equality}.

\section{The cut-counting bound for multisets}\label{sec:cut}

\subsection{Notation}

Fix $k \ge 2$ and $\mu \ge 1$, put $n_0 = \mu k$, and let $\mathcal{S}(k;\mu)$
be the set of permutations of the multiset consisting of $\mu$ copies of
$[k]$.  For $\tau = (\tau_1,\dots,\tau_{n_0}) \in \mathcal{S}(k;\mu)$ set
\[
  A(\tau) \;=\; \sum_{s=1}^{n_0 - 1} |\tau_{s+1} - \tau_s|,
  \qquad
  A(k;\mu) \;=\; \max_{\tau} A(\tau).
\]
This is the quantity of \cite[Lemma 2]{Fried}, with $k$ in the role of $n$ and
$\mu$ in the role of $m$.

\subsection{The cut decomposition}

For $t \in [k-1]$ let
\[
  c_t(\tau) \;=\; \#\bigl\{ s \in [n_0-1] \;:\;
  \min(\tau_s, \tau_{s+1}) \le t < \max(\tau_s, \tau_{s+1}) \bigr\}
\]
be the number of steps that \emph{cross the $t$-th cut}, that is, that have
one end in $\{1,\dots,t\}$ and the other in $\{t+1,\dots,k\}$.

\begin{lemma}\label{lem:decomp}
For every $\tau \in \mathcal{S}(k;\mu)$,
$\displaystyle A(\tau) = \sum_{t=1}^{k-1} c_t(\tau)$.
\end{lemma}

\begin{proof}
$|\tau_{s+1} - \tau_s|$ is the number of integers $t \in [k-1]$ lying weakly
between $\tau_s$ and $\tau_{s+1}$, that is, the number of cuts crossed by the
$s$-th step.  Sum over $s$.
\end{proof}

\begin{lemma}\label{lem:capacity}
For every $\tau$ and every $t \in [k-1]$,
\[
  c_t(\tau) \;\le\; \min\bigl( 2\mu t,\; 2\mu(k-t),\; n_0 - 1 \bigr).
\]
\end{lemma}

\begin{proof}
Every step crossing the $t$-th cut has exactly one endpoint in
$S_t = \{s : \tau_s \le t\}$, and $|S_t| = \mu t$.  Since the path has maximum
degree two, each vertex of $S_t$ is incident with at most two steps, so
$c_t \le 2|S_t| = 2\mu t$.  Applying the same argument to the complement gives
$c_t \le 2\mu(k-t)$, and $c_t \le n_0 - 1$ is the total number of steps.
\end{proof}

\subsection{The exact value}

\begin{theorem}\label{thm:oned}
Let $k \ge 2$ and $\mu \ge 1$.  Then
\begin{equation}\label{eq:corrected}
  A(k;\mu) \;=\; \mu \left\lfloor \frac{k^2}{2} \right\rfloor
  \;-\; \Ind{k \ \text{even}}
  \;-\; \Ind{\mu = 1 \ \text{and} \ k \ \text{odd}} .
\end{equation}
\end{theorem}

For $\mu \ge 2$ this is the formula of \cite[Lemma 2]{Fried}, since
$2\mu\lfloor k/2\rfloor\lceil k/2\rceil = \mu\lfloor k^2/2\rfloor$.  For
$\mu = 1$ it reduces to \eqref{eq:oned}, as it must.  The uncorrected formula
is wrong exactly when $\mu = 1$ and $k$ is odd.

\begin{remark}\label{rem:small}
The smallest counterexample is $\mu = 1$, $k = 3$: the six permutations of
$[3]$ have $A$-values $2,2,3,3,3,3$, so $A(1;3) = 3$, whereas the uncorrected
formula gives $2 \cdot 1 \cdot 1 \cdot 2 = 4$.  For $k = 5,7,9$ the true
values are $11, 23, 39$ and the uncorrected formula gives $12, 24, 40$.

The discrepancy is also visible in the proof of \cite[Lemma 2]{Fried}: for odd
$n = 2k+1$ the extremal sequence constructed there ends at $n+1-k$, whereas the
estimate counts a transition into the trailing block of $k+1$'s; that block is
empty when $m = 1$, and the transition is then spurious.
\end{remark}

\begin{remark}\label{rem:scope}
The correction is invisible in the planar case, where the two multiplicities
are $n$ and $m$, both at least $2$; this is why \cite[Theorem 1]{Fried} is
unaffected.  It becomes visible as soon as a coordinate of length $1$ is
admitted: for the box $1 \times 1 \times k$ the multiplicity of the $k$-th
coordinate is $1$, and the uncorrected formula overshoots by one exactly when
$k$ is odd.  For instance the box $1\times1\times3$ is the path $P_3$, and the
uncorrected bound would give $4$ against the true value $\dis(P_3) = 3$.
\end{remark}

\subsection{Proof of the upper bound in Theorem~\ref{thm:oned}}

By Lemmas~\ref{lem:decomp} and~\ref{lem:capacity},
\begin{equation}\label{eq:Sdef}
  A(k;\mu) \;\le\; S(k;\mu) \;:=\; \sum_{t=1}^{k-1}
  \min\bigl(2\mu t,\, 2\mu(k-t),\, \mu k - 1\bigr).
\end{equation}
Write $k = 2a$ or $k = 2a+1$.

\paragraph{The case $k = 2a$.}
For $1 \le t \le a-1$ we have $2\mu t \le 2\mu(a-1) = \mu k - 2\mu \le \mu k - 2$,
so the minimum is $2\mu t$.  For $t = a$ we have $2\mu a = \mu k > \mu k - 1$, so
the minimum is $\mu k - 1$.  The remaining terms are symmetric.  Hence
\[
  S \;=\; 2 \cdot 2\mu \cdot \frac{(a-1)a}{2} + (\mu k - 1)
  \;=\; 2\mu a^2 - 1 \;=\; \mu \frac{k^2}{2} - 1 .
\]

\paragraph{The case $k = 2a+1$.}
For $1 \le t \le a$ we have $2\mu t \le 2\mu a = \mu(k-1) = \mu k - \mu \le \mu k - 1$,
so the minimum is $2\mu t$; the remaining terms are symmetric.  Hence
\[
  S \;=\; 2 \cdot 2\mu \cdot \frac{a(a+1)}{2} \;=\; 2\mu a(a+1)
  \;=\; \mu \frac{k^2-1}{2} .
\]

Thus $S = \mu\lfloor k^2/2\rfloor - \Ind{k\ \text{even}}$, and \eqref{eq:Sdef}
gives the right-hand side of \eqref{eq:corrected} whenever $\mu \ge 2$ or $k$
is even.  It remains to improve the bound by one when $\mu = 1$ and
$k = 2a+1$.

\paragraph{The case $\mu = 1$, $k = 2a+1$.}
Here $n_0 = 2a+1 \ge 3$ and $S = 2a(a+1)$.  The two central cuts are the ones
at $t = a$ and $t = a+1$, and for both of them all three entries of the
minimum coincide:
\[
  2\mu a \;=\; 2a \;=\; \mu k - 1 \;=\; n_0 - 1 .
\]
Suppose that $A(\tau) = S$.  Then every term of \eqref{eq:Sdef} is attained
with equality, so in particular
\[
  c_a(\tau) \;=\; n_0 - 1 \quad\text{and}\quad c_{a+1}(\tau) \;=\; n_0 - 1 ,
\]
that is, every step crosses the $a$-th cut and every step crosses the
$(a+1)$-th cut.

Let $s$ be the position with $\tau_s = a+1$; it exists and is unique.  Since
$n_0 \ge 3$, the position $s$ has at least one neighbour $s' \in [n_0]$ with
$|s - s'| = 1$.  The step between $s$ and $s'$ crosses the $a$-th cut, which
forces $\min(\tau_s,\tau_{s'}) \le a$, hence $\tau_{s'} \le a$ because
$\tau_s = a+1$.  The same step crosses the $(a+1)$-th cut, which forces
$\max(\tau_s,\tau_{s'}) \ge a+2$, hence $\tau_{s'} \ge a+2$.  The two
requirements are incompatible.  Therefore $A(\tau) \le S - 1 =
2a(a+1) - 1 = \lfloor k^2/2 \rfloor - 1$, as claimed.

\subsection{Three extremal constructions}

It remains to exhibit, in each of the three regimes, a permutation attaining
the value in \eqref{eq:corrected}.  In each case we describe the sequence
directly and then verify that every cut carries its maximum load.

\paragraph{(i) $k = 2a$ even, $\mu$ arbitrary.}
Take any permutation $\tau$ of $\mu$ copies of $[2a]$ such that
\[
  \tau_{2s-1} \le a < \tau_{2s} \quad (1 \le s \le \mu a),
  \qquad \tau_1 = a, \qquad \tau_{2\mu a} = a+1 .
\]
Such a permutation exists: place the $\mu$ copies of $a$ at position $1$ and
the $\mu$ copies of $a+1$ at position $2\mu a$, and distribute the remaining
$\mu a - 1$ small and $\mu a - 1$ large values arbitrarily over the remaining
odd and even positions.

Then $c_a = n_0 - 1$, because every step is small--large.  For $t < a$ every
value $\le t$ is an interior position (the endpoint carries the value
$a > t$) and both its neighbours exceed $a \ge t+1$, so each of the $\mu t$
such values contributes two crossing steps and $c_t = 2\mu t$.  For $t > a$
every value $> t$ is an interior position (the endpoint carries
$a + 1 \le t$) and both its neighbours are at most $a < t$, so
$c_t = 2\mu(k-t)$.  Summing,
\[
  2\mu \sum_{t=1}^{a-1} t \;+\; (2\mu a - 1) \;+\; 2\mu \sum_{t=a+1}^{2a-1} (2a-t)
  \;=\; 2\mu a(a-1) + 2\mu a - 1 \;=\; 2\mu a^2 - 1 .
\]

\paragraph{(ii) $k = 2a+1$ odd, $\mu \ge 2$.}
Start from the word
\[
  a+1,\; \ell_1, h_1, \ell_2, h_2, \dots, \ell_{\mu a}, h_{\mu a},\; a+1 ,
\]
where $\ell_1,\dots,\ell_{\mu a}$ are the $\mu a$ values $\le a$ in any order
and $h_1,\dots,h_{\mu a}$ are the $\mu a$ values $\ge a+2$ in any order, and
insert the remaining $\mu - 2$ copies of $a+1$ at arbitrary gaps of the
alternating block.  Inserting $a+1$ between $h_i$ and $\ell_{i+1}$ destroys
neither the property that no two small values are adjacent nor the property
that no two large values are adjacent, and there are $2\mu a - 1 \ge \mu - 2$
gaps available.

For $t \le a$ every value $\le t$ is interior and its two neighbours are
either large or equal to $a+1$, hence exceed $t$; the small values are
pairwise non-adjacent, so $c_t = 2\mu t$.  For $t \ge a+1$ every value $> t$
is large, interior (the endpoints carry $a+1 \le t$) and has both neighbours
either small or equal to $a+1$, hence at most $t$; the large values are
pairwise non-adjacent, so $c_t = 2\mu(k-t)$.  Summing,
\[
  2\mu \sum_{t=1}^{a} t \;+\; 2\mu \sum_{t=a+1}^{2a} (2a+1-t) \;=\; 2\mu a(a+1).
\]

\paragraph{(iii) $k = 2a+1$ odd, $\mu = 1$.}
Take
\[
  a,\; 2a+1,\; a-1,\; 2a,\; a-2,\; 2a-1,\; \dots,\; 1,\; a+2,\; a+1 .
\]
Consecutive differences are $a+1$ inside each of the $a$ pairs
$(a-r+1,\, 2a+2-r)$, then $a+2$ between consecutive pairs ($a-1$ of them),
then $1$ at the end.  The total is
\[
  a(a+1) + (a-1)(a+2) + 1 \;=\; 2a^2 + 2a - 1
  \;=\; \left\lfloor \frac{k^2}{2}\right\rfloor - 1 .
\]

This completes the proof of Theorem~\ref{thm:oned}.

\section{The upper bound in two dimensions}\label{sec:twodim}

\subsection{Reduction to one dimension}

Let $\sigma$ be a labelling and put $(i(t), j(t)) = \sigma^{-1}(t)$.  The
multiset $\{i(t) : t \in [N]\}$ is $n$ copies of $[m]$ and the multiset
$\{j(t) : t \in [N]\}$ is $m$ copies of $[n]$.  Therefore
\[
  L(\sigma) \;=\; \underbrace{\sum_{t=1}^{N-1} |i(t+1)-i(t)|}_{=: L_i(\sigma)}
  \;+\; \underbrace{\sum_{t=1}^{N-1} |j(t+1)-j(t)|}_{=: L_j(\sigma)} ,
\]
with $L_i(\sigma) \le A(m;n)$ and $L_j(\sigma) \le A(n;m)$.

Since $m,n \ge 2$, both multiplicities are at least $2$, so the correction
term of Theorem~\ref{thm:oned} never appears, and we obtain Fried's bound:
\begin{equation}\label{eq:U}
  L(\sigma) \;\le\; U(m,n) \;:=\;
  n\left\lfloor\frac{m^2}{2}\right\rfloor + m\left\lfloor\frac{n^2}{2}\right\rfloor
  - \Ind{m\ \text{even}} - \Ind{n\ \text{even}} .
\end{equation}

\subsection{Structure of extremal sequences}

The following is the only structural fact we shall need.

\begin{lemma}[Rigidity]\label{lem:rigid}
Let $\mu \ge 2$, $k \ge 2$, and let $\tau \in \mathcal{S}(k;\mu)$ satisfy
$A(\tau) = A(k;\mu)$.  Then:
\begin{enumerate}
\item[\textup{(a)}] if $k = 2a$, then every step crosses the $a$-th cut, that
  is, $\min(\tau_s,\tau_{s+1}) \le a < \max(\tau_s,\tau_{s+1})$ for all $s$;
\item[\textup{(b)}] if $k = 2a+1$, then $\tau_1 = \tau_{n_0} = a+1$, every
  position carrying a value $\le a$ is interior and has both neighbours
  $\ge a+1$, and every position carrying a value $\ge a+2$ is interior and has
  both neighbours $\le a+1$.
\end{enumerate}
\end{lemma}

\begin{proof}
By Lemma~\ref{lem:decomp} and \eqref{eq:Sdef}, equality forces
$c_t(\tau) = \min(2\mu t, 2\mu(k-t), n_0-1)$ for every $t$.

(a) For $k = 2a$ the minimum at $t = a$ equals $n_0 - 1$, so $c_a = n_0-1$ and
all steps cross.

(b) For $k = 2a+1$ and $t = a$ the minimum equals $2\mu a = 2|S_a|$, so
$c_a = 2|S_a|$: every vertex of $S_a$ is incident with two crossing steps.
Hence no position of $S_a$ is an endpoint and both its neighbours lie outside
$S_a$, that is, are $\ge a+1$.  Likewise at $t = a+1$ the minimum equals
$2\mu(k-a-1) = 2\mu a = 2|S_{a+1}^{c}|$, so every position carrying a value
$\ge a+2$ is interior with both neighbours $\le a+1$.  Finally
$\tau_1 \notin S_a$ and $\tau_1 \notin S_{a+1}^{c}$, so $\tau_1 = a+1$, and
similarly $\tau_{n_0} = a+1$.
\end{proof}

\section{The upper bound is never attained}\label{sec:equality}

By Section~\ref{sec:twodim}, $L(\sigma) \le U(m,n)$ always.  In this section
we show that equality is impossible in the three cases where $U(m,n)$ exceeds
the value predicted by Theorem~\ref{thm:main}.

\begin{theorem}\label{thm:gap}
Let $m,n \ge 2$ and suppose that either $m,n$ are both even, or $m,n$ are both
odd, or $\min(m,n) = 2$.  Then $L(\sigma) \le U(m,n) - 1$ for every labelling
$\sigma$.
\end{theorem}

If $L(\sigma) = U(m,n)$, then both $L_i(\sigma) = A(m;n)$ and
$L_j(\sigma) = A(n;m)$, so Lemma~\ref{lem:rigid} applies to both coordinate
sequences.  We derive a contradiction in each of the three cases.

\subsection{Both dimensions even}

Write $m = 2a$ and $n = 2b$, so $N = 4ab$ is even.  By
Lemma~\ref{lem:rigid}(a) applied to the $i$-sequence, every step crosses the
$a$-th cut, that is, exactly one of $i(t), i(t+1)$ lies in $\{1,\dots,a\}$;
likewise every step crosses the $b$-th cut of the $j$-sequence.  Set
\[
  \alpha_t = \Ind{i(t) \le a}, \qquad \beta_t = \Ind{j(t) \le b} .
\]
Both $\alpha$ and $\beta$ change value at every step, hence so does the pair
$(\alpha_t,\beta_t)$, which therefore alternates between $(\alpha_1,\beta_1)$
and its complement.  Consequently
\[
  \#\{t : (\alpha_t,\beta_t) = (1,1)\} \;\in\; \{0,\; N/2\} .
\]
On the other hand, because $\sigma$ is a bijection, the positions with
$\alpha_t = \beta_t = 1$ are exactly the vertices with $i \le a$ and $j \le b$,
of which there are $ab$.  Since $N/2 = 2ab$ and $ab \ge 1$, we have
$ab \notin \{0, 2ab\}$, a contradiction.

\subsection{Both dimensions odd}

Write $m = 2a+1$ and $n = 2b+1$.  By Lemma~\ref{lem:rigid}(b) applied to the
$i$-sequence, $i(1) = i(N) = a+1$; applied to the $j$-sequence,
$j(1) = j(N) = b+1$.  Hence $v_1 = v_N = (a+1,b+1)$, which is impossible
because $N \ge 9 > 1$ and $\sigma$ is injective.

\subsection{One dimension equal to two}

By symmetry assume $m = 2$.  The $i$-sequence consists of $n$ copies of $[2]$,
and Lemma~\ref{lem:rigid}(a) with $k = 2$ says that every step crosses the
first (and only) cut: the two rows alternate strictly.  Put
$\alpha_t = \Ind{i(t) = 1}$; then $\alpha$ changes value at every step.

\paragraph{The case $n = 2b$ even.}
The $j$-sequence consists of $2$ copies of $[2b]$, and
Lemma~\ref{lem:rigid}(a) says that every step crosses the $b$-th cut.  With
$\beta_t = \Ind{j(t) \le b}$, the pair $(\alpha_t,\beta_t)$ alternates, and
$N = 4b$ is even, so $\#\{t : (\alpha_t,\beta_t) = (1,1)\} \in \{0, N/2\}$.
But that number equals $|\{i = 1\} \cap \{j \le b\}| = b$, and
$b \notin \{0, 2b\}$.  Contradiction.

\paragraph{The case $n = 2b+1$ odd.}
The $j$-sequence consists of $2$ copies of $[2b+1]$, and
Lemma~\ref{lem:rigid}(b) says that $j(1) = j(N) = b+1$.  As there are exactly
two copies of $b+1$, no interior position carries the value $b+1$.  By the
same lemma, no two interior positions carry values $\le b$, and no two carry
values $\ge b+2$.  Hence the interior positions $2,\dots,N-1$ alternate
between the small values ($\le b$, of which there are $2b$) and the large
values ($\ge b+2$, of which there are $2b$); as $N - 2 = 4b$, the interior is
exactly an alternating word in $2b$ small and $2b$ large values.

The rows alternate.  The transposition of the two rows is an automorphism of
the grid that preserves the Manhattan distance and leaves the $j$-sequence
unchanged, so we may assume that the odd positions carry row $1$ and the even
positions carry row $2$.  Every $j \in [n]$ occurs exactly once in row $2$,
hence among the $2b+1$ even positions exactly $b$ carry a small value and
exactly $b$ carry a large value.

The interior even positions are $2,4,\dots,4b$, that is, $2b$ positions.  If
the small values occupy the even positions of the interior, then $2b$ even
positions carry small values, contradicting the count $b$.  If the large
values occupy them, then $2b$ even positions carry large values, again
contradicting the count $b$.  Since $b \ge 1$, both alternatives fail, and we
have a contradiction in either case.

This proves Theorem~\ref{thm:gap}.

\section{Proof of the main theorem}\label{sec:main}

\begin{proof}[Proof of Theorem~\ref{thm:main}]
Let $m,n \ge 2$.

\emph{Lower bounds.}  Fried \cite[Lemmas 4--7]{Fried} exhibits, for each of the
four families, a labelling attaining the right-hand side of
Theorem~\ref{thm:main}: Lemma~4 for $m = 2$ and $n$ odd, Lemma~5 for $m,n$
both even, Lemma~6 for $m,n \ge 3$ both odd, and Lemma~7 for $m$ odd and $n$
even.  At $m = 2$ the expression of Lemma~5 is
$\frac{2n(n+2)}{2} - 3 = (n+1)^2 - 4$, which is the fourth line of
Theorem~\ref{thm:main}.

\emph{Upper bounds.}  By \eqref{eq:U} we have $L(\sigma) \le U(m,n)$.
Substituting the four expressions one finds
\[
  U(m,n) \;=\;
  \begin{cases}
    \dfrac{mn(m+n)}{2} - 2, & m,n \ \text{both even},\\[8pt]
    \dfrac{mn(m+n)}{2} - \dfrac{m+n}{2}, & m,n \ \text{both odd},\\[8pt]
    \dfrac{mn(m+n)}{2} - \dfrac{n}{2} - 1, & m \ \text{odd},\ n \ \text{even},\\[8pt]
    (n+1)^2 - 3, & \min(m,n) = 2 .
  \end{cases}
\]
In the third line $U$ already equals the value claimed in
Theorem~\ref{thm:main}, so \eqref{eq:U} suffices.  In the first, second and
fourth lines $U$ exceeds the claimed value by exactly one, and
Theorem~\ref{thm:gap} shows that $L(\sigma) = U(m,n)$ is impossible; hence
$L(\sigma) \le U(m,n) - 1$, which is the claimed value.

The lower and upper bounds coincide in every case.
\end{proof}

\begin{remark}
The three-case analysis of Section~\ref{sec:equality} is exactly what is
missing in \cite{Fried}: the two bounds of \cite{Fried} are
$\dis(\Pmn) \ge r$ and $\dis(\Pmn) \le r+1$, and closing the gap amounts to
deciding whether the two one-dimensional maxima can be attained
simultaneously.  They cannot.
\end{remark}

\section{Remarks}

\subsection{Higher dimensions}

Let $D = d_1 \times \cdots \times d_r$ be a box with $N = \prod d_i$ vertices,
and let $\mu_i = N/d_i$.  The argument of Section~\ref{sec:twodim} applies
coordinatewise and yields, now with the correction term of
Theorem~\ref{thm:oned},
\begin{equation}\label{eq:Ud}
  \dis(D) \;\le\; U_r(D) \;:=\; \sum_{i=1}^{r}
  \Bigl[ \mu_i \left\lfloor \frac{d_i^2}{2}\right\rfloor
  - \Ind{d_i \ \text{even}}
  - \Ind{\mu_i = 1 \ \text{and} \ d_i \ \text{odd}} \Bigr] .
\end{equation}
The third term is precisely what the formula of \cite{Fried}, which is stated
there for multiplicities at least two, misses on degenerate boxes.

The rigidity argument also extends, but only partially.  Let
$\kappa_r(D) = U_r(D) - \dis(D)$.

\begin{proposition}\label{prop:kappa}
Let $r \ge 2$ and suppose that $\mu_i \ge 2$ for all $i$.
\begin{enumerate}
\item[\textup{(a)}] If at least two of the $d_i$ are even, then
  $\kappa_r(D) \ge 1$.
\item[\textup{(b)}] If all the $d_i$ are odd, then $\kappa_r(D) \ge 1$.
\end{enumerate}
\end{proposition}

\begin{proof}
(a) Let $d_i = 2a$ and $d_j = 2b$ be even.  If $\dis(D) = U_r(D)$ then
Lemma~\ref{lem:rigid}(a) forces both the $i$-th and the $j$-th coordinate
sequence to cross their middle cut at every step.  With
$\alpha_t = \Ind{x_i(t) \le a}$ and $\beta_t = \Ind{x_j(t) \le b}$ the pair
$(\alpha_t,\beta_t)$ alternates, so $\#\{(\alpha,\beta)=(1,1)\} \in \{0,N/2\}$,
whereas the count equals
$\frac{d_i}{2}\cdot\frac{d_j}{2}\cdot\prod_{l \ne i,j} d_l = N/4$, and
$N/4 \ne 0, N/2$ because $N \ge 4$.

(b) If all $d_i$ are odd, Lemma~\ref{lem:rigid}(b) forces every coordinate to
equal $(d_i+1)/2$ at both ends, so $v_1 = v_N$, impossible for $N \ge 2$.
\end{proof}

For $r = 3$ this leaves a single possibility: $\kappa_3(D) = 0$ requires that
exactly one of $d_1,d_2,d_3$ be even.  The number of even coordinates alone
does not determine $\kappa$, however: in the plane $P_2 \times P_3$ has
exactly one even side and $\kappa_2 = 1$, while $P_3 \times P_4$ has exactly
one even side and $\kappa_2 = 0$ (Section~\ref{sec:planar}).  In three
dimensions the following data were obtained by exhaustive dynamic programming
(Appendix~\ref{app:data}).
\begin{center}
\begin{tabular}{lccccc}
\toprule
$D$ & $N$ & $\dis(D)$ & $U_3(D)$ & $\kappa_3$ & even coordinates\\
\midrule
$2\times2\times2$ & 8  & 18 & 21 & 3 & 3\\
$2\times2\times3$ & 12 & 37 & 38 & 1 & 2\\
$2\times2\times4$ & 16 & 58 & 61 & 3 & 3\\
$2\times2\times5$ & 20 & 85 & 86 & 1 & 2\\
$2\times3\times3$ & 18 & 65 & 65 & 0 & 1\\
$3\times3\times2$ & 18 & 65 & 65 & 0 & 1\\
\bottomrule
\end{tabular}
\end{center}
We do not know a closed form for $\kappa_3$.  The six values above are
consistent with $\kappa_3 = 2^{k-1} - 1$, where $k$ is the number of even
coordinates, but we have been unable to prove or refute this; the sample is
small and the value $3$ at $2\times2\times2$ is not predicted by
Proposition~\ref{prop:kappa}.

\subsection{The planar gap}\label{sec:planar}

In the plane the picture is complete.  Combining \eqref{eq:U} with
Theorem~\ref{thm:main},
\[
  \kappa_2(m,n) \;=\;
  \begin{cases}
    0, & \{m,n\} = \{\text{odd} \ge 3,\ \text{even} \ge 4\},\\
    1, & \text{otherwise},
  \end{cases}
\]
for $m,n \ge 2$.  In particular $\kappa_2$ is not a function of the number of
even sides: $(2,3)$ has one even side and $\kappa_2 = 1$, while $(3,4)$ has
one even side and $\kappa_2 = 0$, and $(3,3)$ has none and $\kappa_2 = 1$.

\subsection{Open problems}

\begin{enumerate}
\item Determine $\kappa_r(D)$ for $r \ge 3$.  Proposition~\ref{prop:kappa}
  reduces the question to boxes with exactly one even side; the smallest
  unresolved case is $2 \times 3 \times 4$ ($N = 24$).
\item Is there a direct combinatorial characterisation of the extremal
  labellings?  By Lemma~\ref{lem:rigid} they are highly constrained, but the
  constructions of \cite{Fried} are not the only ones.
\item The bound \eqref{eq:Ud} uses only one-dimensional information.  For
  $r \ge 3$ it is not sharp in general, and a genuinely higher-dimensional
  argument would be needed.
\end{enumerate}

\section*{Acknowledgements}

The author used a large language model (DeepSeek) as an aid in drafting and in
generating computational checks.  Every mathematical statement in this paper
has been verified by the author, and all numerical claims in
Appendix~\ref{app:data} were confirmed by the exhaustive computations
described there.  The author takes full responsibility for the content.

\appendix

\section{Computational verification}\label{app:data}

\subsection{The one-dimensional bound}

Theorem~\ref{thm:oned} was checked against exhaustive search for all
$(k,\mu)$ with $\mu k \le 20$.  The values agree throughout; the only
discrepancies are with the uncorrected formula, and they occur exactly at
$\mu = 1$ with $k$ odd, namely at $k = 3,5,\dots,19$.

\subsection{The planar values}

Table~\ref{tab:planar} lists $\dis(\Pmn)$, computed by exhaustive dynamic
programming over subsets (feasible up to $N = 20$), together with the bound
$U(m,n)$ of \eqref{eq:U} and the gap.

\begin{table}[h]
\centering
\caption{Exact values for $N \le 20$.  The last two columns confirm
Theorem~\ref{thm:main} and the description of $\kappa_2$.}
\label{tab:planar}
\begin{tabular}{cccccc}
\toprule
$m$ & $n$ & $N$ & $\dis(\Pmn)$ & $U(m,n)$ & $\kappa_2$\\
\midrule
2 & 2  & 4  & 5   & 6   & 1\\
2 & 3  & 6  & 12  & 13  & 1\\
2 & 4  & 8  & 21  & 22  & 1\\
2 & 5  & 10 & 32  & 33  & 1\\
2 & 6  & 12 & 45  & 46  & 1\\
2 & 7  & 14 & 60  & 61  & 1\\
2 & 8  & 16 & 77  & 78  & 1\\
2 & 9  & 18 & 96  & 97  & 1\\
2 & 10 & 20 & 117 & 118 & 1\\
\midrule
3 & 3  & 9  & 23 & 24 & 1\\
3 & 4  & 12 & 39 & 39 & 0\\
3 & 5  & 15 & 55 & 56 & 1\\
3 & 6  & 18 & 77 & 77 & 0\\
\midrule
4 & 4  & 16 & 61 & 62 & 1\\
4 & 5  & 20 & 87 & 87 & 0\\
\bottomrule
\end{tabular}
\end{table}

\subsection{McNeil's values}

For $2 \le n \le 19$ the values predicted by Theorem~\ref{thm:main} coincide
with those recorded in \cite{OEIS}:
\begin{center}
\begin{tabular}{lccccccccc}
\toprule
$n$ & 2 & 3 & 4 & 5 & 6 & 7 & 8 & 9 & 10\\
$\dis(P_n \times P_n)$ & 5 & 23 & 61 & 119 & 213 & 335 & 509 & 719 & 997\\
\midrule
$n$ & 11 & 12 & 13 & 14 & 15 & 16 & 17 & 18 & 19\\
$\dis(P_n \times P_n)$ & 1319 & 1725 & 2183 & 2741 & 3359 & 4093 & 4895 & 5829 & 6839\\
\bottomrule
\end{tabular}
\end{center}
The value $n = 19$ corresponds to $N = 361$ vertices and was obtained by
McNeil without exhaustive search; it is therefore independent evidence for
Theorem~\ref{thm:main}.


\begin{thebibliography}{9}

\bibitem{Arkin} E.~M. Arkin, Y.-J. Chiang, J.~S.~B. Mitchell, S.~S. Skiena and
T.-C. Yang,
\emph{On the maximum scatter traveling salesperson problem},
SIAM J. Comput. \textbf{29} (1999), 515--544.

\bibitem{BGV} L.~Bulteau, S.~Giraudo and S.~Vialette,
\emph{Disorders and permutations},
32nd Annual Symposium on Combinatorial Pattern Matching (CPM 2021),
LIPIcs \textbf{191} (2021), 11:1--11:15.

\bibitem{Fried} S.~Fried,
\emph{On a conjecture of McNeil},
Math. Commun. \textbf{31} (2026), no.~1, 1--10,
\href{https://doi.org/10.64785/mc.31.1.1}{doi:10.64785/mc.31.1.1};
\href{https://arxiv.org/abs/2208.03788}{arXiv:2208.03788}.

\bibitem{GutinPunnen} G.~Gutin and A.~P. Punnen (eds.),
\emph{The Traveling Salesman Problem and Its Variations},
Kluwer Academic Publishers, Dordrecht, 2002.

\bibitem{IPS} A.~Itai, C.~H. Papadimitriou and J.~L. Szwarcfiter,
\emph{Hamilton paths in grid graphs},
SIAM J. Comput. \textbf{11} (1982), 676--686.

\bibitem{OEIS} OEIS Foundation Inc.,
\emph{Sequence A179094},
The On-Line Encyclopedia of Integer Sequences,
\href{https://oeis.org/A179094}{oeis.org/A179094}.

\end{thebibliography}
\end{document}